\documentclass{amsart}
\usepackage{amsmath}
\usepackage{amssymb}
\usepackage{graphicx}
\usepackage{color}
\usepackage{enumerate}
\usepackage[english,vietnamese]{babel}
\newtheorem{theorem}{Theorem}[section]

\newtheorem{conjecture}{Conjecture}
\newtheorem{lemma}[theorem]{Lemma}
\newtheorem{corollary}[theorem]{Corollary}
\newtheorem{proposition}[theorem]{Proposition}
\theoremstyle{definition}

\newtheorem{remark}{Remark}
\newenvironment{Proof}{{\textit{Proof}.}\ }{~$\square$\vspace{0.2truecm}}

\newcommand{\Aut}{\mbox{\rm Aut}}

\newcommand{\Z}{\mathbb{Z}}
\newcommand{\C}{\mathbb{C}}

\newcommand{\M}{{\rm M}}
\newcommand{\I}{{\rm I}}

\newcommand{\GL}{{\rm GL}}
\newcommand{\SL}{{\rm SL}}
\newcommand{\Gal}{{\rm Gal}}

\begin{document}

\title[Solvable maximal subgroups]{A note on solvable maximal subgroups\\ in subnormal subgroups of  $\GL_n(D)$}
\author[Hu\`{y}nh Vi\d{\^{e}}t Kh\'{a}nh]{Hu\`{y}nh Vi\d{\^{e}}t Kh\'{a}nh}\thanks{The second author was funded by Vietnam National Foundation for Science and Technology Development (NAFOSTED) under Grant No. 101.04-2016.18}	
\address{Faculty of Mathematics and Computer Science, VNUHCM - University of Science, 227 Nguyen Van Cu Str., Dist. 5, Ho Chi Minh City, Vietnam.} 
\email{huynhvietkhanh@gmail.com; bxhai@hcmus.edu.vn} 
\author[B\`{u}i Xu\^{a}n H\h{a}i]{B\`{u}i Xu\^{a}n H\h{a}i}	
\maketitle
\selectlanguage{english}

\begin{abstract} Let $D$ be a non-commutative division ring, and $G$ be a subnormal subgroup of $\GL_n(D)$. Assume additionally that the center of $D$ contains at least five elements if $n>1$. In this note, we  show that if $G$ contains a non-abelian solvable maximal subgroup, then $n=1$ and $D$ is a cyclic algebra of prime degree over the center.
\end{abstract}

\section{Introduction}

In the theory of skew linear groups, one of unsolved difficult problems is that whether the general skew linear group over a division ring contains maximal subgroups. In \cite{akbri}, the authors conjectured that for $n\geq 2$ and a division ring $D$, the group $\GL_n(D)$ contains no solvable maximal subgroups. In \cite{dorbidi2011}, this conjecture was shown to be true for non-abelian solvable maximal subgroups. In this paper, we consider the following more general conjecture.

\begin{conjecture}\label{conj:1} Let $D$ be a division ring, and $G$ be a non-central subnormal subgroup of $\GL_n(D)$. If $n\geq 2$, then $G$ contains no solvable maximal subgroups.
\end{conjecture}

We note that this conjecture is not true if $n=1$. Indeed, it was proved in \cite{akbri} that the subgroup $\C^*\cup \C^* j$ is a solvable maximal subgroup of the multiplicative group $\mathbb{H}^*$ of the division ring of real quaternions $\mathbb{H}$. In this note, we show that Conjecture \ref{conj:1} is true for non-abelian solvable maximal subgroups of $G$, that is, we prove that $G$ contains no non-abelian solvable maximal subgroups. This fact generalizes the main result in \cite{dorbidi2011} and it is a consequence of Theorem \ref{theorem_3.7} in the text.  

Throughout this note, we denote by $D$ a division ring with center $F$ and by $D^*$ the multiplicative group of $D$. For a positive integer $n$, the symbol $\M_n(D)$ stands for the matrix ring of degree $n$ over $D$. We identify $F$ with $F\I_n$ via the ring isomorphism $a\mapsto a\I_n$, where $\I_n$ is the identity matrix of degree $n$. If $S$ is a subset of $\M_n(D)$, then $F[S]$ denotes the subring of $\M_n(D)$ generated by the set $S\cup F$. Also, if $n=1$, i.e., if $S\subseteq D$, then $F(S)$ is the division subring of $D$ generated by $S\cup F$. Recall that a division ring $D$ is \textit{locally finite} if for every finite subset $S$ of $D$, the division subring $F(S)$ is a finite dimensional vector space over $F$. If $H$ and $K$ are two subgroups in a group $G$, then $N_K(H)$ denotes the set of all elements $k\in K$ such that $k^{-1}Hk\leq H$, i.e., $N_K(H)=K\cap N_G(H)$. If $A$ is a ring or a group, then $Z(A)$ denotes the center of $A$.

Let $V =D^n= \left\{ {\left( {{d_1},{d_2}, \ldots ,{d_n}} \right)\left| {{d_i} \in D} \right.} \right\}$. If $G$ is a subgroup of $\GL_n(D)$, then $V$ may be viewed as $D$-$G$ bimodule. Recall that a subgroup $G$ of $\GL_n(D)$ is  \textit{irreducible} (resp. \textit{reducible, completely reducible}) if $V$ is irreducible (resp. reducible, completely reducible) as $D$-$G$ bimodule. If $F[G]=\M_n(D)$, then   $G$ is  \textit{absolutely irreducible} over $D$. An irreducible subgroup $G$ is \textit{imprimitive} if there exists an integer $ m \geq 2$ such that $V =  \oplus _{i = 1}^m{V_i}$ as left $D$-modules and for any $g \in G$ the mapping $V_i \to V_ig$ is a permutation of the set $\{V_1, \cdots, V_m\}$. If $G$ is irreducible and not imprimitive, then $G$ is \textit{primitive}. 

\section{Auxiliary lemmas}

\begin{lemma}\label{lemm_2.1}
	Let $D$ be a division ring with center $F$, and $M$ be a subgroup of $\GL_n(D)$. If $M/M\cap F^*$ is a locally finite group, then $F[M]$ is a locally finite dimensional vector space over $F$.
\end{lemma}

\begin{Proof} Take any finite subset $\{x_1,x_2,\dots,x_k\}\subseteq F[M]$ and write
	$$x_i=f_{i_1}m_{i_1}+f_{i_2}m_{i_2}+\cdots+f_{i_s}m_{i_s}.$$
	Let $G=\left\langle m_{i_j}:1\leq i\leq k,1\leq j\leq s\right\rangle$ be the subgroup of $M$ generated by all $m_{i_j}$.
	Since $M/M\cap F^*\cong MF^*/F^*$ is locally finite, the group $GF^*/F^*$ is finite. Let $\{y_1,y_2,\dots,y_t\}$ be a transversal of $F^*$ in $GF^*$ and set
	$$R=Fy_1+Fy_2+\cdots+Fy_t.$$ 
	Then, $R$ is a finite dimensional vector space over $F$ containing $\{x_1,x_2,\dots,x_k\}$. 
\end{Proof}

\begin{lemma}\label{lemma_2.2}
	Every locally solvable periodic group is locally finite.
\end{lemma}

\begin{Proof} 
	Let $G$ be a locally solvable periodic group, and $H$ be a finitely generated subgroup of $G$. Then, $H$ is solvable with derived series of length $n\geq1$, say,
	$$1=H^{(n)}\unlhd H^{(n-1)}\unlhd\cdots\unlhd H'\unlhd H.$$
	We shall prove that $H$ is finite by induction on $n$. For if $n=1$, then $H$ is a finitely generated periodic abelian group, so it is finite. Suppose $n>1$. It is clear that $H/H'$ is a finitely generated periodic abelian group, so it is finite. Hence, $H'$ is finitely generated. By induction hypothesis, $H'$ is finite, and as a consequence,  $H$ is finite.
\end{Proof}

\begin{lemma}\label{lemma_2.3}
	Let $D$ be a division ring with center $F$, and $G$ be a subnormal subgroup of $D^*$. If $G$ is solvable-by-finite, then $G\subseteq F$.
\end{lemma}

\begin{Proof} 
	Let $A$ be a solvable normal subgroup of finite index in $G$. Since $G$ is subnormal in $G$, so is $A$. By \cite[14.4.4]{scott}, we have $A\subseteq F$. This implies that $G/Z(G)$ is finite, so $G'$ is finite too \cite[Lemma 1.4, p. 115]{passman_77}. Therefore, $G'$ is a finite subnormal subgroup of $D^*$. In view of \cite[Theorem 8]{her}, it follows that $G'\subseteq F$, hence $G$ is solvable. Again by \cite[14.4.4]{scott}, we conclude that $G\subseteq F$.
\end{Proof}

For our further use, we also need one result of Wehrfritz which will be restated in the following lemma for readers' convenience. 

\begin{lemma}\cite[Proposition 4.1]{wehrfritz_07}\label{lemma_2.4}
	Let $D=E(A)$ be a division ring generated as such by its metabelian subgroup $A$ and its division subring $E$ such that $E\leq C_D(A)$. Set $H=N_{D^*}(A)$, $B=C_A(A')$, $K=E(Z(B))$, $H_1=N_{K^*}(A)=H\cap K^*$, and let $T$ be the maximal periodic normal subgroup of $B$.
	\begin{enumerate}[(1)]
		\item If $T$ has a quaternion subgroup $Q=\left\langle i,j\right\rangle $ of order $8$ with $A=QC_A(Q)$, then $H=Q^+AH_1$, where $Q^+=\left\langle Q,1+j,-(1+i+j+ij)/2\right\rangle$. Also, $Q$ is normal in $Q^+$ and $Q^+/{\left\langle -1,2\right\rangle}\cong\Aut Q\cong Sym(4)$.	
		\item If $T$ is abelian and contains an element $x$ of order $4$ not in the center of $B$, then $H=\left\langle x+1\right\rangle AH_1$.
		\item In all other cases, $H=AH_1$.
	\end{enumerate}
\end{lemma}

\section{Maximal subgroups in subnormal subgroups of $\GL_n(D)$}

\begin{proposition}\label{proposition_3.1}
	Let $D$ be a division ring with center $F$, and $G$ be a subnormal subgroup of $D^*$. If $M$ is a non-abelian solvable-by-finite maximal subgroup of $G$, then $M$ is abelian-by-finite and $[D:F]<\infty$.
\end{proposition}

\begin{Proof} 
	Since $M$ is maximal in $G$ and $M\subseteq F(M)^*\cap G\subseteq G$, either $M = F(M)^*\cap G$ or $G \subseteq F(M)^*$. The first case implies that $M$ is a solvable-by-finite subnormal subgroup of $F(M)^*$, which yields $M$ is abelian by Lemma \ref{lemma_2.3}, a contradiction. Therefore, the second case must occur, i.e., $G \subseteq F(M)^*$. By Stuth's theorem (see e.g. \cite[14.3.8]{scott}), we conclude that $F(M)=D$. Let $N$ be a solvable normal subgroup of finite index in $M$. First, we assume that $N$ is abelian, so $M$ is abelian-by-finite. In view of \cite[Corollary 24]{wehrfritz_89}, the ring $F[N]$ is a Goldie ring, and hence it is an Ore domain whose skew field of fractions  coincides with $F(N)$. Consequently, any $\alpha\in F(N)$ may be written in the form $\alpha=pq^{-1},$ where $q, p\in F[N]$ and $q\ne0$. The normality of $N$ in $M$ implies that $F[N]$ is normalized by $M$. Thus, for any $m\in M$, we have
		$$m\alpha m^{-1}=mpq^{-1}m^{-1}=(mpm^{-1})(m^{-1}qm)^{-1}\in F(N).$$
	In other words, $L:=F(N)$ is a subfield of $D$ normalized by $M$. Let $\{x_1,x_2,\ldots,x_k\}$ be a transversal of $N$ in $M$ and set
	$$\Delta=Lx_1+Lx_2+\cdots+Lx_k.$$
	Then, $\Delta$ is a domain with $\dim_L\Delta\leq k$, so $\Delta$ is a division ring that is finite dimensional over its center. It is clear that $\Delta$ contains $F$ and $M$, so $D=\Delta$ and $[D:F]<\infty$. 
	
	Next, we suppose that $N$ is a non-abelian solvable group with derived series of length $s\geq1$. Then we have such a series
		$$1=N^{(s)}\unlhd N^{(s-1)}\unlhd N^{(s-2)}\unlhd\cdots\unlhd N'\unlhd  N \unlhd M.$$
	If we set $A=N^{(s-2)}$, then $A$ is a non-abelian metabelian normal subgroup of $M$. By the same arguments as above, we conclude that $F(A)$ is normalized by $M$ and we have $M\subseteq N_G(F(A)^*)\subseteq G$. By the maximality of $M$ in $G$, either $N_G(F(A)^*)=M$ or $N_G(F(A)^*)=G$. If the first case occurs, then $G\cap F(A)^*$ is a subnormal subgroup of $F(A)^*$ contained in $M$. Since $M$ is solvable-by-finite, so is $G\cap F(A)^*$. By Lemma \ref{lemma_2.3}, $A\subseteq G\cap F(A)^*$ is abelian, a contradiction. We may therefore assume that $N_G(F(A))=G$, which says that $F(A)$ is normalized by $G$. In view of Stuth's theorem, we have $F(A)=D$. From this we conclude that $Z(A)=F^*\cap A$ and $F=C_D(A)$. Set $H=N_{D^*}(A)$,  $B=C_A(A')$, $K=F(Z(B))$, $H_1=H\cap K^*$, and $T$ to be the maximal periodic normal subgroup of $B$. Then $H_1$ is an abelian group, and $T$ is a characteristic subgroup of $B$ and hence of $A$. In view of Lemma \ref{lemma_2.4}, we have three possible cases:
	
	\bigskip
	
	\textit{Case 1:} $T$ is not abelian. 
	
	\bigskip
	
	Since $T$ is normal in $M$, we conclude that $M\subseteq N_G(F(T)^*)\subseteq G$. By the maximality of $M$ in $G$, either $M = N_G(F(T)^*)$ or $G= N_G(F(T)^*)$. The first case implies that $F(T)^*\cap G$ is subnormal in $F(T)^*$ contained in $M$. Again by Lemma~ \ref{lemma_2.3}, it follows that $T\subseteq F(T)\cap G$ is abelian, a contradiction. Thus, we may assume that $G= N_G(F(T)^*)$, which implies that $F(T)=D$ by Stuth's theorem. By Lemma~ \ref{lemma_2.2}, $T$ is locally finite. In view of Lemma \ref{lemm_2.1}, we conclude that $D=F(T)=F[T]$ is a locally finite division ring. Since  $M$ is solvable-by-finite, it contains no non-cyclic free subgroups. In view of \cite[Theorem 3.1]{hai-khanh}, it follows $[D:F]<\infty$ and $M$ is abelian-by-finite.
	
	\bigskip	
	
	\textit{Case 2:} $T$ is abelian and contains an element $x$ of order $4$ not in the center of $B=C_A(A')$.
	
	\bigskip
	
	It is clear that $x$ is not contained in $F$. Because $x$ is of finite order, the field $F(x)$ is algebraic over $F$.  Since $\left\langle x\right\rangle$ is a $2$-primary component of $T$, it is a characteristic subgroup of $T$ (see the proof of \cite[Theorem 1.1, p. 132]{wehrfritz_07}). Consequently, $\left\langle x\right\rangle$ is a normal subgroup of $M$. Thus, all elements of the set $x^M:=\{m^{-1}xm\vert m\in M\}\subseteq F(x)$ have the same minimal polynomial over $F$. This implies $|x^M|<\infty$, so $x$ is an $FC$-element, and consequently, $[M:C_M(x)]<\infty$. Setting $C=Core_M(C_M(x))$, then $C\unlhd M$ and $[M:C]$ is finite. Since $M$ normalizes $F(C)$, we have $M\subseteq N_G(F(C)^*) \subseteq G$. By the maximality of $M$ in $G$, either $N_G(F(C)^*)=M$ or $N_G(F(C)^*)=G$. The last case implies that $F(C)=D$, and consequently, $x\in F$, a contradiction. Thus, we may assume that $N_G(F(C)^*)=M$. From this, we conclude that $G\cap F(C)^*$ is a subnormal subgroup of $F(C)^*$ which is contained in $M$. Thus, $C\subseteq G\cap F(C)^*$ is contained in the center of $F(C)$ by \cite[14.4.4]{scott}. Therefore, $C$ is an abelian normal subgroup of finite index in $M$. By the same arguments used in the first paragraph we conclude that $[D:F]<\infty$.
	
	\bigskip	
	
	\textit{Case 3:} $H=AH_1$.
	
	\bigskip
	
	Since $A'\subseteq H_1\cap A$, we have $H/H_1\cong A/A\cap H_1$ is abelian, and hence $H'\subseteq H_1$. Since $H_1$ is abelian, $H'$ is abelian too. Moreover, $M\subseteq H$, it follows that $M'$ is also abelian. In other words, $M$ is a metabelian group, and the conclusions follow from \cite[Theorem 3.3]{hai-tu}.
\end{Proof}

Let $D$ be a division ring, and $G$ be a subnormal subgroup of $D^*$. It was showed in \cite[Theorem 3.3]{hai-tu} that if $G$ contains a non-abelian metabelian maximal subgroup, then $D$ is cyclic of prime degree. The following theorem generalizes this phenomenon.

\begin{theorem} \label{theorem_3.2}
	Let $D$ be a division ring with center $F$, and $G$ be a subnormal subgroup of $D^*$. If $M$ is a non-abelian solvable maximal subgroup of $G$, then the following conditions hold:
	\begin{enumerate}[(i)]
		\item There exists a maximal subfield $K$ of $D$ such that $K/F$ is a finite Galois extension with $\mathrm{Gal}(K/F)\cong M/K^*\cap G\cong \mathbb{Z}_p$ for  some prime $p$, and $[D:F]=p^2$. 
		\item The subgroup $K^*\cap G$ is the $FC$-center. Also, $K^*\cap G$ is  the Fitting subgroup of $M$. Furthermore, for any $x\in M\setminus K$, we have $x^p\in F$ and $D=F[M]=\bigoplus_{i=1}^pKx^i$.
	\end{enumerate}
\end{theorem}

\begin{Proof}
	By Proposition \ref{proposition_3.1}, it follows that $[D:F]<\infty$. Since $M$ is solvable, it contains no non-cyclic free subgroups. In view of \cite[Theorem 3.4]{hai-tu}, we have $F[M]=D$, there exists a maximal subfield $K$ of $D$  containing $F$ such that $K/F$ is a Galois extension, $N_G(K^*)=M $, $K^*\cap G$ is the Fitting normal subgroup of $M$ and it is the $FC$-center, and $M/K^*\cap G\cong\Gal(K/F)$ is a finite simple group of order $[K:F]$. Since $M/K^*\cap G$ is solvable and simple, one has $M/K^*\cap G\cong\Gal(K/F)\cong \Z_p$, for some prime number $p$. Therefore, $[K:F]=p$ and $[D:F]=p^2$. For any $x\in M\backslash K$, if $x^p\not\in F$, then by the fact that $F[M]=D$, we conclude that $C_M(x^p)\ne M$. Moreover, since $x^p\in K^*\cap G$, it follows that $\left\langle x,K^*\cap G\right\rangle \leq C_M(x^p)$. In other words, $C_M(x^p)$ is a subgroup of $M$ strictly containing $K^*\cap G$. Because $M/K^*\cap G$ is simple, we have $C_M(x^p)= M$, a contradiction. Therefore $x^p\in F$. Furthermore, since $x^p\in K$ and $[D:K]_r=p$, we conclude $D=\bigoplus_{i=1}^{p-1}Kx_i$.
\end{Proof}

 Also, the authors in \cite{nassab14} showed that if $D$ is an infinite division ring, then $D^*$ contains no polycyclic-by-finite maximal subgroups. In the following corollary, we will see that every subnormal subgroup of $D^*$ does not contain non-abelian polycyclic-by-finite maximal subgroups.

\begin{corollary}
	Let $D$ be a division ring with center $F$, $G$ be a subnormal subgroup of $D^*$, and $M$ be a non-abelian maximal subgroup of $G$. Then $M$ cannot be  finitely generated solvable-by-finite. In particular, $M$ cannot be polycyclic-by-finite.
\end{corollary}

\begin{Proof}
	Suppose that $M$ is solvable-by-finite. Then by Proposition \ref{proposition_3.1} , we conclude that $[D:F]<\infty$. In view of \cite[Corollary 3]{mah2000}, it follows that $M$ is not finitely generated. The rest of the corollary is clear.
\end{Proof}

\begin{theorem} \label{theorem_3.4} 
	Let $D$ be a non-commutative locally finite division ring with center $F$, and $G$ be a subnormal subgroup of $\GL_n(D)$, $n\geq 1$. If $M$ is a non-abelian solvable maximal subgroup of $G$, then $n=1$ and all conclusions of Theorem \ref{theorem_3.2} hold.
\end{theorem}

\begin{Proof}
	By \cite[Theorem 3.1]{hai-khanh}, there exists a maximal subfield $K$ of $\M_n(D)$  containing $F$ such that $K^*\cap G$ is a normal subgroup of $M$ and $M/K^*\cap G$ is a finite simple group of order $[K:F]$. Since $M/K^*\cap G$ is solvable and simple, we conclude $M/K^*\cap G\cong \Z_p$, for some prime number $p$. It follows that $[K:F]=p$ and $[M_n(D):F]=p^2$, from which we have $n=1$. Finally, all conclusions follow from Theorem \ref{theorem_3.2}.
\end{Proof}

\begin{lemma}\label{lemma_3.5}
	Let $R$ be a ring, and $G$ be a subgroup of $R^*$. Assume that $F$ is a central subfield of $R$ and that $A$ is a maximal abelian subgroup of $G$ such that $K=F[A]$ is normalized by $G$. Then $F[G]=\oplus _{g \in T}{Kg}$ for every transversal $T$ of $A$ in $G$.
\end{lemma}

\begin{Proof} 
	For the proof of this lemma, we use the similar techniques as in the proof of  \cite[Lemma 3.1]{dorbidi2011}.
	Since $K$ is normalized by $G$, it follows that $F[G]=\sum\nolimits_{g \in T}{Kg}$ for every transversal $T$ of $A$ in $G$. Therefore, it suffices to prove that every finite subset $\{g_1,g_2,\dots,g_n\}\subseteq T$ is  linearly independent over $K$. Assume by contrary that there exists such a non-trivial relation 
		$$k_1g_1+k_2g_2+\cdots+k_ng_n=0.$$
	Clearly, we can suppose that all $k_i$ are non-zero, and that $n$ is minimal. If $n=1$, then there is nothing to prove, so we can suppose $n>1$. Since the cosets $Ag_1$ and $Ag_2$ are disjoint, we have $g_1^{-1}g_2\not\in A=C_G(A)$. So, there exists an element $x\in A$ such that $g_1^{-1}g_2x\ne xg_1^{-1}g_2$. For each $1\leq i\leq n$, if we set $x_i=g_ixg_i^{-1}$, then $x_1\ne x_2$. Since $G$ normalizes $K$, it follows $x_i\in K$ for all $1\leq i\leq n$. Now, we have 
		$$(k_1g_1+\cdots+k_ng_n)x-x_1(k_1g_1+\cdots+k_ng_n)= 0.$$
	By definition, we have $x_ig_i=g_ix$, and $x$, $x_i\in K$ for all $i$. Recall that $K=F[A]$ is commutative, so from the last equality
		$$\left( {{k_1}{g_1}x + {k_2}{g_2}x +  \cdots  + {k_n}{g_n}x} \right) - \left( {{x_1}{k_1}{g_1} + {x_1}{k_2}{g_2} +  \cdots  + {x_1}{k_n}{g_n}} \right) = 0,$$
	it follows
		$$\left( {{k_1}{x_1}{g_1} + {k_2}{x_2}{g_2} +  \cdots  + {k_n}{x_n}{g_n}} \right) - \left( {{k_1}{x_1}{g_1} + {k_2}{x_1}{g_2} +  \cdots  + {k_n}{x_1}{g_n}} \right) = 0.$$
	Consequently, we get
		$$\left( {{x_2} - {x_1}} \right){k_2}{g_2} +  \cdots  + \left( {{x_n} - {x_1}} \right){k_n}{g_n} = 0,$$
	which is a non-trivial relation with less than $n$ summands because $x_1\ne x_2$, a contradiction. Therefore, $T$ is linearly independent over $K$.	
\end{Proof}

\begin{remark}\label{rem.subnormal}
	In view of \cite[Theorem 11]{mah98}, if $D$ is a division ring with at least five elements and $n\geq 2$, then any non-central subnormal subgroup of $\GL_n(D)$ contains $\SL_n(D)$ and hence is normal.
\end{remark}
 
\begin{theorem} \label{theorem_3.6}
	Let $D$ be non-commutative division ring with center $F$, and $G$ be a subnormal subgroup of $\GL_n(D)$, $n\geq2$. Assume additionally that $F$ contains at least five elements. If $M$ is a solvable maximal subgroup of $G$, then $M$ is abelian.
\end{theorem}

\begin{Proof} 
	If $G\subseteq F$, then there is nothing to prove. Thus, we may assume that $G$ is non-central, hence $\SL_n(D)\subseteq G$ and $G$ is normal in $\GL_n(D)$ by Remark \ref{rem.subnormal}. Setting $R=F[M]$, then $M\subseteq R^*\cap G \subseteq G$. By the maximality of $M$ in $G$, either $R^*\cap G=M$ or $G\subseteq R^*$.  We need to consider two possible cases:
	\bigskip	
	
	\textit{Case 1:} $R^*\cap G=M$.
	
	\bigskip
	
	The normality of $G$ implies that $M$ is a normal subgroup of $R^*$. If $M$ is reducible, then by \cite[Lemma 1]{kiani_07}, it contains a copy of $D^*$. It follows that $D^*$ is solvable, and hence it is commutative, a contradiction. We may therefore assume that $M$ is irreducible. Then $R$ is a prime ring by \cite[1.1.14]{shirvani-wehrfritz}. So, in view of \cite[Theorem 2]{lanski_81}, either $M\subseteq Z(R)$ or $R$ is a domain. If the first case occurs, then we are done. Now, suppose that $R$ is a domain. By \cite[Corollary 24]{wehrfritz_89}, we conclude that $R$ is a Goldie ring, and thus $R$ is an Ore domain. Let $\Delta_1$ be the skew field of fractions of $R$, which is contained in $\M_n(D)$ by \cite[5.7.8]{shirvani-wehrfritz}. Since $M\subseteq \Delta_1^*\cap G\subseteq G$, either $G\subseteq\Delta_1^*$ or $M=\Delta_1^*\cap G$. The first case occurs implies that $\Delta_1$ contains $F[\SL_n(D)]$. Thus, if $G\subseteq\Delta_1^*$, then by the Cartan-Brauer-Hua Theorem for the matrix ring (see e.g. \cite[Theorem D]{dorbidi2011}), one has $\Delta_1=\M_n(D)$, which is impossible since $n\geq 2$. Thus the second case must occur, i.e., $M=\Delta_1\cap G$, which yields $M$ is normal in $\Delta_1^*$. Since $M$ is solvable, it is contained in $Z(\Delta_1)$ by \cite[14.4.4]{scott}, so $M$ is abelian.
	
	\bigskip	
	
	\textit{Case 2:} $G\subseteq R^*$.
	
	\bigskip
	
	In this case, Remark \ref{rem.subnormal} yields  $\SL_n(D)\subseteq R^*$. Thus, by the Cartan-Brauer-Hua Theorem for the matrix ring, one has $R=F[M]=\M_n(D)$. It follows by \cite[Theorem A]{wehrfritz86} that $M$ is abelian-by-locally finite. Let $A$ be a maximal abelian normal subgroup of $M$ such that $M/A$ is locally finite. Then \cite[1.2.12]{shirvani-wehrfritz} says that $F[A]$ is a semisimple artinian ring. The Wedderburn-Artin Theorem implies that
	$$ F[A] \cong \M_{n_1}(D_1)\times \M_{n_2}(D_2)\cdots\times \M_{n_s}(D_s),$$
	where $D_i$ are division $F$-algebras, $1\leq i\leq s$. Since $F[A]$ is abelian, $n_i = 1$ and $K_i:=D_i=Z(D_i)$ are fields that contain $F$ for all $i$. Therefore, 
	$$F[A]\cong K_1\times K_2\cdots\times K_s.$$
	If $M$ is imprimitive, then by \cite[Lemma 2.6]{hai-khanh}, we conclude that $M$ contains a copy of $\SL_r(D)$ for some $r\geq 1$. This fact cannot occur: if $r>1$, then $\SL_r(D)$ is unsolvable; if $r=1$, then $D'$ is solvable and hence $D$ is commutative, a contradiction. Thus, $M$ is primitive, and \cite[Proposition 3.3]{dorbidi2011} implies that $F[A]$ is an integral domain, so $s=1$. It follows that $K:=F[A]$ is a subfield of $\M_n(D)$ containing $F$. Again by \cite[Proposition 3.3]{dorbidi2011}, we conclude that $$L:=C_{\M_n(D)}(K)=C_{\M_n(D)}(A)\cong \M_m(\Delta_2)$$ for some division $F$-algebra $\Delta_2$. Since $M$ normalizes $K$, it also normalizes $L$. Therefore, we have $M\subseteq N_G(L^*)\subseteq G$. By the maximality of $M$ in $G$, either $M= N_G(L^*)$ or $G=N_G(L^*)$. The last case implies that $L^*\cap G$ is normal in $\GL_n(D)$. By Remark \ref{rem.subnormal}, either $L^*\cap G\subseteq F$ or $\SL_n(D)\subseteq L^*\cap G$. If the first case occurs, then $A\subseteq F$ because $A$ is contained in $L^*\cap G$. If the second case occurs, then by the Cartan-Brauer-Hua Theorem for the matrix, one has $L=\M_n(D)$. It follows that $K=F[A]\subseteq F$, which also implies that $A\subseteq F$. Thus, in both case we have $A\subseteq F$. Consequently, $M/M\cap F^*$ is locally finite, and hence $D$ is a locally finite division ring by Lemma \ref{lemm_2.1}. If $M$ is non-abelian, then by Theorem \ref{theorem_3.4}, we conclude that $n=1$, a contradiction. Therefore $M$ is abelian in this case. Now, we consider the case $M= N_G(L^*)$, from which we have $L^*\cap G \subseteq M$. In other words, $L^*\cap G$ is a solvable normal subgroup of $\GL_m(\Delta_2)$. From this, we may conclude that $L^*\cap G$ is abelian:  if $m>1$, then in view of Remark \ref{rem.subnormal}, one has $L^*\cap G\subseteq Z(\Delta_2)$ or $\SL_m(\Delta_2)\subseteq L^*\cap G$, but the latter cannot happen since $\SL_m(\Delta_2)$ is unsolvable; if $m=1$ then $L=\Delta_2$, and according to \cite[14.4.4]{scott}, we conclude that $L^*\cap G\subseteq Z(\Delta_2)$. In short, we have $L^*\cap G$ is an abelian normal subgroup of $M$ and $M/L^*\cap G$ is locally finite. By the maximality of $A$ in $M$, one has $A=L^*\cap G$. Because we are in the case $L^*\cap G \subseteq M$, it follows that $L^*\cap G=L^*\cap M$. Consequently,  $$A=L^*\cap M=C_{\GL_n(D)}(A)\cap M =C_M(A),$$ which means $A$ is a maximal abelian subgroup of $M$.
	
	By Lemma \ref{lemma_3.5}, $F[M]=\oplus _{m \in T}{Km}$ for some transversal $T$ of $A$ in $M$. Thus, for any $x\in L$, there exist $k_1,k_2,\dots,k_t\in K$ and $m_1,m_2,\dots,m_t\in T$ such that $x=k_1m_1+k_2m_2+\cdots+k_tm_t$. Take an arbitrary element $a\in A$, since $xa=ax$, it follows that
	$$k_1m_1a+k_2m_2a+\cdots+k_tm_ta=ak_1m_1+ak_2m_2+\cdots+ak_tm_t.$$ 
	By the normality of $A$ in $M$, there exist $a_i\in A$ such that $m_ia=a_im_i$ for all $1\leq i\leq t$. Moreover, we have $a$ and $a_i$'s are in $K$ which is a field, the equality implies
	$$k_1a_1m_1+k_2a_2m_2+\cdots+k_ta_tm_t=k_1am_1+k_2am_2+\cdots+k_tam_t,$$ from which it follows that
	$$ (k_1a_1-k_1a)m_1+(k_2a_2-k_2a)m_2+\cdots+(k_ta_t-k_ta)m_t=0.$$ 
	Since $\{m_1,m_2,\dots,m_t\}$ is linearly independent over $K$, one has $a=a_1=\cdots=a_t$. Consequently, $m_ia=am_i$ for all $a\in A$, and thus $m_i\in C_M(A)=A$ for all $1\leq i\leq t$. This means $x\in K$, and hence $L=K$ and $K$ is a maximal subfield of $\M_n(D)$.
	
	Next, we prove that $M/A$ is simple. Suppose that $N$ is an arbitrary normal subgroup of $M$ properly containing $A$. Note that by the maximality of $A$ in $M$, we conclude that $N$ is non-abelian. We claim that $Q:=F[N]=\M_n(D)$. Indeed, since $N$ is normal in $M$, we have $M\subseteq N_G(Q^*)\subseteq G$, and hence either $N_G(Q^*)=M$ or $N_G(Q^*)=G$. First, we suppose the former case occurs. Then $Q^*\cap G\subseteq M$, hence $Q^*\cap G$ is a solvable normal subgroup of $Q^*$. In view of \cite[Proposition 3.3]{dorbidi2011}, $Q$ is a prime ring. It follows by \cite[Theorem 2]{lanski_81} that either $Q^*\cap G\subseteq Z(Q)$ or $Q$ is a domain. If the first case occurs, then $N\subseteq Q^*\cap G$ is abelian, which contradicts to the choice of $N$. If $Q$ is a domain, then by Goldie's theorem, it is an Ore domain. Let $\Delta_3$ be the skew field of fractions of $Q$, which is contained in $\M_n(D)$ by \cite[5.7.8]{shirvani-wehrfritz}. Because $M$ normalizes $Q$, it also normalizes $\Delta_3$, from which we have $M\subseteq N_G(\Delta_3^*)\subseteq G$. Again by the maximality of $M$ in $G$, either $N_G(\Delta_3^*)=M$ or $N_G(\Delta_3^*)=G$. The first case implies that $\Delta_3^*\cap G$ is a solvable normal subgroup of $\Delta_3^*$. Consequently, $N\subseteq\Delta_3^*\cap G$ is abelian by \cite[14.4.4]{scott}, a contradiction. If $N_G(\Delta_3^*)=G$, then $\Delta_3=\M_n(D)$ by the Cartan-Brauer-Hua Theorem for the matrix ring, which is impossible since $n\geq2$. Therefore, the case $N_G(Q^*)=M$ cannot occur. Next, we consider the case $N_G(Q^*)=G$. In this case we have $Q^*\cap G \unlhd \GL_n(D)$, and hence either $Q^*\cap G\subseteq F$ or $\SL_n(D)\subseteq Q^*\cap G$ by Remark \ref{rem.subnormal}. The first case cannot occur since $Q^*\cap G$ contains $N$, which is non-abelian. Therefore, we have $\SL_n(D)\subseteq Q^*$. By the Cartan-Brauer-Hua theorem for the matrix ring, we conclude $Q=\M_n(D)$ as claimed. In other words, we have $F[N]=F[M]=\M_n(D)$.
	
	For any $m\in M\subseteq F[N]$, there exist $f_1,f_2,\dots,f_s\in F$ and $n_1,n_2,\dots,n_s\in N$ such that 
	$$m=f_1n_1+f_2n_2+\cdots+f_sn_s.$$
	Let $H=\left\langle n_1,n_2,\dots,n_s \right\rangle $ be the subgroup of $N$ generated by $n_1,n_2,\dots,n_s$. Set $B=AH$ and $S=F[B]$. Recall that $A$ is a maximal abelian subgroup of $M$. Thus, if $B$ is abelian, then $A=B$ and hence $H\subseteq A$. Consequently, $m\in F[A]=K$, from which it follows that $m\in K^*\cap M=A\subseteq N$. Now, assume that $B$ is non-abelian, and we will prove that $m$ also belongs to $N$ in this case. Since $M/A$ is locally finite, $B/A$ is finite. Let $\{x_1,\ldots,x_k\}$ be a transversal of $A$ in $B$. The maximality of $A$ in $M$ implies that $A$ is a maximal abelian subgroup of $B$, and that $A$ is also normal in $B$. By Lemma \ref{lemma_3.5},
	$$S=Kx_1\oplus Kx_2\oplus\cdots\oplus Kx_k,$$ 
	which says that $S$ is an artinian ring. Since $C_{\M_n(D)}(A)=C_{\M_n(D)}(K)=L$ is a field, in view of \cite[Proposition 3.3]{dorbidi2011}, we conclude that $A$ is irreducible. Because $B$ contains $A$, by definition, it is irreducible too. It follows by \cite[1.1.14]{shirvani-wehrfritz} that $S$ is a prime ring. Now, $S$ is both prime and artinian, so it is simple and $S\cong \M_{n_0}(\Delta_0)$ for some division $F$-algebra $\Delta_0$. If we set $F_0=Z(\Delta_0)$, then $Z(S)=F_0$. Since $B$ is abelian-by-finite, the group ring $FB$ is a PI-ring by \cite[Lemma 11, p. 176]{passman_77}. Thus, as a hommomorphic image of $FB$, the ring $S=F[B]$ is also a PI-ring. By Kaplansky's theorem (\cite[Theorem 3.4, p. 193]{passman_77}), we conclude that $[S:F_0]<\infty$. Since $K$ is a maximal subfield of $\M_n(D)$, it is also maximal in $S$. From this, we conclude that $F_0\subseteq C_S(K)= K$, and that $K$ is a finite extension field over $F_0$.

	Recall that $A$ is normal in $B$, so for any $b\in B$,  the mapping $\theta_b:K\to K$ given by $\theta_b(x)=bxb^{-1}$ is well defined. It is clear that $\theta_b$ is an $F_0$-automorphism of $K$. Thus, the mapping $$\psi:B\to \Gal(K/F_0)$$ defined by $\psi(b)=\theta_b$ is a group homomorphism with $$\mathrm{ker}\psi=C_{S^*}(K^*)\cap B=C_{S^*}(A)\cap B=C_B(A)=A.$$ Since $F_0[B]=S$, it follows that $C_S(B)=F_0$. Therefore, the fixed field of $\psi(B)$ is $F_0$, and hence $K/F_0$ is a Galois extension. By the fundamental theorem of Galois theory, one has $\psi$ is a surjective homomorphism. Hence, $ B/A\cong \Gal(K/F_0)$.
	
	Setting $M_0=M\cap S^*$, then $m\in M_0$, $B\subseteq M_0$, and $F_0[M_0]=F_0[B]=S$. The conditions $F_0\subseteq K$ and $F\subseteq F_0$ implies that $K=F[A]=F_0[A]$. It is clear that $A$ is a maximal abelian subgroup of $M_0$, and that $A$ is also normal in $M_0$. If $M_0/A$ is infinite, then there exists an infinite transversal $T$ of $A$ in $M_0$ such that $S=F_0[M_0]=\oplus_{m\in T}Km$ by Lemma \ref{lemma_3.5}. It follows that $[S:K]=\infty$, a contradiction. Therefore, $M_0/A$ must be finite.  Replacing $B$ by $M_0$ in the preceding paragraph, we also conclude that $ M_0/A\cong \Gal(K/F_0)$. Consequently, $B/A\cong\Gal(K/F_0)\cong M_0/A$. The conditions $B\subseteq M_0$ and  $B/A \cong M_0/A$ imply $B=M_0$. Hence, $m\in M_0= B\subseteq N$.  Since $m$ was chosen arbitrarily, it follows that $M=N$, which implies the simplicity of $M/A$. Since $M/A$ is simple and solvable, one has $M/A\cong \Z_p$, for some prime number $p$. By Lemma \ref{lemma_3.5}, it follows $\dim_K\M_n(D)=|M/A|=p$, which forces $n=1$, a contradiction. 
\end{Proof}

Now, we are ready to get the main result of this note which gives in particular, the positive answer to Conjecture \ref{conj:1} for non-abelian case.
\begin{theorem} \label{theorem_3.7}
	Let $D$ be a non-commutative division ring with center $F$, $G$ a subnormal subgroup of $\GL_n(D)$. Assume additionally that $F$ contains at least five elements if $n>1$. If $M$ is a non-abelian solvable maximal subgroup of $G$, then $n=1$ and the following conditions hold:
	\begin{enumerate}[(i)]
		\item There exists a maximal subfield $K$ of $D$ such that $K/F$ is a finite Galois extension with $\mathrm{Gal}(K/F)\cong M/K^*\cap G\cong \mathbb{Z}_p$ for  some prime $p$, and $[D:F]=p^2$. 
		\item The subgroup $K^*\cap G$ is the $FC$-center. Also, $K^*\cap G$ is  the Fitting subgroup of $M$. Furthermore, for any $x\in M\setminus K$, we have $x^p\in F$ and $D=F[M]=\bigoplus_{i=1}^pKx^i$.
	\end{enumerate}
\end{theorem}

\begin{Proof}
 Combining Theorem \ref{theorem_3.2} and Theorem \ref{theorem_3.6}, we get the result.
\end{Proof}


\begin{thebibliography}{}
	
	\bibitem{akbri} S. Akbari, R. Ebrahimian, H. Momenaei Kermani, and A. Salehi Golsefidy. Maximal subgroups of $\GL_n(D)$. \textit{J. Algebra} \textbf{259} (2003), 201-225.
	
	\bibitem{dorbidi2011}  H. R. Dorbidi, R. Fallah-Moghaddam, and M. Mahdavi-Hezavehi. Soluble maximal subgroups in $\GL_n(D)$. \textit{J. Algebra Appl.} (6) \textbf{10} (2011), 1371-1382.  
	
	\bibitem{draxl}  P. Draxl. \textit{Skew Fields}. London Math. Soc. Lecture Note Ser. \textbf{81} (Cambridge University Press, 1983).
	
	\bibitem{hai-tu}  B. X. Hai and N. A. Tu. On multiplicative subgroups in division rings. \textit{J. Algebra  Appl.} (3) \textbf{15} (2016), 1650050 (16 pages).
	
	\bibitem{hai-khanh}  B. X. Hai and H. V. Khanh. On free subgroups in maximal subgroups of skew linear groups. \textit{arXiv:1808.08453v1 [math.RA].}
	
	\bibitem{her} I. N. Herstein. Multiplicative commutators in division rings. \textit{Israel J. Math.} \textbf{31} (1978), 180-188.
	
	\bibitem{kiani_07} D. Kiani. Polynomial identities and maximal subgroups of skew linear groups. \textit{Manuscripta Math.} \textbf{124} (2007), 269-274.
	
	\bibitem{lanski_81} C. Lanski. Solvable subgroups in prime rings.  \textit{Proc. Amer. Math. Soc.} \textbf{82} (1981), 533-537. 
	
	\bibitem{mah98} M. Mahdavi-Hezavehi and S. Akbari. Some special subgroups of $\GL_n(D)$. \textit{Algebra Colloq.} (4) \textbf{5} (1998), 361-370.
	
	\bibitem{mah2000} M. Mahdavi-Hezavehi, M. G. Mahmudi, and S. Yasamin. Finitely generated subnormal subgroups of $\GL_n(D)$ are central. \textit{J. Algebra} \textbf{225} (2000), 517-521.
	
	\bibitem{nassab14} M. Ramezan-Nassab and D. Kiani. Nilpotent and polycyclic-by-finite maximal subgroups of skew linear groups. \textit{J. Algebra} \textbf{399} (2014), 269-276.
		
	\bibitem{passman_77} D. S. Passman. \textit{The algebraic structure of group rings}. (New York: Wiley- Interscience Publication, 1977).
	
	\bibitem{scott} R. W. Scott. \textit{Group Theory}. (Dover Publication, INC, 1987).
	
	\bibitem{shirvani-wehrfritz} M. Shirvani and B. A. F. Wehrfritz. \textit{Skew Linear Groups}. (Cambridge Univ. Press, 1986).
	
	\bibitem{wehrfritz86}  B. A. F. Wehrfritz. Soluble normal subgroups of skew linear groups. \textit{J. Pure Appl. Algebra} \textbf{42} (1986), 95-107.
	
	\bibitem{wehrfritz_07} B. A. Wehrfritz. Normalizers of nilpotent subgroups of division rings. \textit{Q. J. Math.} \textbf{58} (2007), 127-135.
	
	\bibitem{wehrfritz_89} B. A. F. Wehrfritz. Goldie subrings of Artinian rings generated by groups. \textit{Q. J. Math. Oxford} \textbf{40} (1989), 501-512.
\end{thebibliography}
\end{document}